\documentclass[twoside,12pt]{article}
\usepackage{graphicx,amscd,amsfonts,amsmath,amsthm,amssymb,latexsym,amsfonts,color}
\usepackage{epsfig,float,epstopdf,array,bm,cite,cases,array,multirow,relsize,graphicx,float}
\usepackage[bookmarksnumbered,colorlinks,plainpages]{hyperref}
\usepackage[table]{xcolor}
\newtheorem{theorem}{\bf Theorem}
\newtheorem{lemma}{\bf Lemma}

\def\h{\hspace{-0.2cm}}

\begin{document}
\title{\bf A new block alternating splitting iteration method for solving a class of  two-by-two block complex linear systems}
\author{\small\bf Davod Khojasteh Salkuyeh$^{\dag\ddag}$\thanks {\noindent Corresponding author.}, Maeddeh Pourbagher$^{\dag}$  \\[2mm]
	\textit{{\small $^{\dag}$Faculty of Mathematical Sciences, University of Guilan, Rasht, Iran}} \\
	\textit{{\small $^{\ddag}$Center of Excellence for Mathematical Modelling, Optimization and Combinational}}\\
	\textit{{\small Computing (MMOCC), University of Guilan, Rasht, Iran}}\\
	\textit{{\small Emails: khojasteh@guilan.ac.ir,mapourbagher@gmail.com}} 
}
\date{}
\maketitle
\vspace{-0.3cm}

\noindent\hrulefill\\
{\bf Abstract.}  A  block alternating splitting  iteration (BASI) method is presented for solving the system arising from  the finite element discretization of the distributed optimal control problem with time-periodic parabolic equations. We prove that the BASI	 method is unconditionally convergent.  We derive the BASI preconditioner and  present an estimation formula for the parameter of the BASI preconditioner. Numerical results are presented to verify  the efficiency of both the BASI  method and the BASI preconditioner.  Comparison with some existing methods are also given. 
 \\[-3mm]


\noindent\hrulefill

\pagestyle{myheadings}\markboth{D.K. Salkuyeh, M. Pourbagher}{Modified block alternating splitting}
\thispagestyle{empty}

\section{Introduction} \label{Sec1}\rm
We are concerned  with the following distributed control problem
\begin{eqnarray*}
	\min_{y,u}  \dfrac{1}{2} \int_{0}^{T} \int_{\Omega} | y(x,t)- y_{d}(x,t)|^{2} dx dt&\h+\h& \dfrac{\nu}{2} \int_{0}^{T} \int_{\Omega} |u(x,t)|^{2}dxdt,\\
	\text{s.t:}\quad\dfrac{\partial}{\partial t}y(x,t)- \triangle y(x,t) &\h=\h& u(x,t) \qquad \text{in} \quad  Q_{T},\\
	y(x,t) &\h=\h& 0      \qquad \qquad \text{on} \quad \Sigma_{T},\\
	y(x,0)&\h=\h& y(x,T) \quad \quad \text{on} \quad \partial \Omega,\\
	u(x,0)&\h=\h&u(x,T) \quad \quad \text{in} \quad  \Omega,
\end{eqnarray*}
where $\Omega$ is an open and bounded domain in $\mathbb{R}^{d}$ for $d \in \{1,2,3\}$ with Lipschitz-continuous boundary $\partial \Omega$.
The space-time cylinder and its lateral surface for $T>0$ are difined as  $Q_{T}=\Omega \times (0,T)$ and $\Sigma_{T}=\partial \Omega \times (0,T)$, respectively. Here $y_{d}(x,t)$ is a given target and $\nu >0$ is a cost or regularization parameter. We assume that the function $y_{d}(x,t)$ is time-harmonic, which means that 
$$y_{d}(x,t)=y_{d}(x) e^{i \omega t}\quad \text{with} \quad\omega= \frac{2 m \pi}{T},$$ 
for some $m \in \mathbb{Z}$. It can be proved that there is a time-periodic solution for 
the original problem as 
\[
y(x,t)=y(x) e^{i \omega t}\quad \text{and} \quad u(x,t)=u(x) e^{i \omega t},
\]
where $y(x)$ and $u(x)$ solve the following time-independent optimal control problem:
\begin{eqnarray*}
	\min_{y,u} \dfrac{1}{2} \int_{\Omega} | y(x)- y_{d}(x)|^{2} dx &\h+\h& \dfrac{\nu}{2} \int_{\Omega} |u(x)|^{2}dx,\\
	\text{s.t:}\quad 
	i \omega y(x)- \triangle y(x) &\h=\h& u(x,t) \qquad \text{in} \quad \Omega,\\
	y(x) &\h=\h& 0      \qquad \qquad \text{on} \quad \partial \Omega.
\end{eqnarray*}
Substituting $y_{d}(x,t)$, $y(x,t)$ and $u(x,t)$ in the model problem gives an optimal control problem.  Discretizing the obtained optimal control problem  using the finite element method  results in an optimization problem. If we define the corresponding Lagrangian functional with the Lagrange multiplier $\bar{p}$, then the first order necessary conditions which are also sufficient for the existence of a solution is obtained as the following  
complex system of linear equations   
\begin{equation}\label{MAIN}
	Ax\equiv \begin{pmatrix}
		M & \sqrt{\nu} (K-i \omega M)\\ 
		\sqrt{\nu} (K+i \omega M) & -M\\ 
	\end{pmatrix} \begin{pmatrix}
		\bar{y}\\ 
		\bar{q}\\ 
	\end{pmatrix} =\begin{pmatrix}
		M \bar{y}_{d}\\ 
		0\\ 
	\end{pmatrix}\equiv b, 
\end{equation}
with $\bar{u}={\bar{p}}/{\nu}$ and $\bar{q}={\bar{p}}/{\sqrt{\nu}}$, where $\bar{y}$, $\bar{y}_d$, and $\bar{u}$ are the coefficient vectors of the corresponding finite element basis functions of $y(x,t)$, $y_{d}(x,t)$ and $u(x,t)$, respectively. 
Moreover,  $M, K \in \mathbb{R}^{m \times m}$ are the mass and stiffness matrices, respectively. Besides, the linear system (\ref{MAIN}) can be transformed to 
\begin{equation}\label{Axel}
	\hat{A} x \equiv \begin{pmatrix}
		M & -\sqrt{\nu} (K-i \omega M)\\ 
		\sqrt{\nu} (K+i \omega M) & M\\ 
	\end{pmatrix} \begin{pmatrix}
		\bar{y}\\ 
		\hat{q}\\ 
	\end{pmatrix} =\begin{pmatrix}
		M \bar{y}_{d}\\ 
		0\\ 
	\end{pmatrix}\equiv b, 
\end{equation}
and $\hat{q}=-{\bar{p}}/{\sqrt{\nu}}$.
Since the coefficient matrix of \eqref{MAIN} is of large size, iterative methods such as the ones based on the Krylov subspace are employed to solve the system. In order to accelerate the convergence of iterative methods suitable preconditioners should be manipulated.  
In \cite{Krendl}, Krendl et al. proposed the real block diagonal and the alternative indefinite preconditioners 
for the system (\ref{MAIN}).  The block alternating splitting (BAS) iteration method was presented by Zheng et al. in \cite{Zheng1} for solving the system (\ref{MAIN}) which can be written as
\begin{equation}\label{BAS}
	\left\{ \begin{array}{l}
		\left( \alpha V +H_{1}\right) x^{(k+\frac{1} {2})}=\left( \alpha V - S_{1} \right)  x^{(k)}+P_{1}b, \\
		\left( \alpha V+H_{2} \right) x^{(k+1)}=\left(  \alpha V -S_{2} \right)  x^{(k+\frac{1}{2})}+P_{2} b,
	\end{array}\right.
\end{equation}
where $\alpha>0$, $V \in \mathbb{R}^{2m \times 2m}$ is a symmetric positive definite (SPD) matrix,
\begin{eqnarray*}
	H_1= \begin{pmatrix}
		M& 0 \\ 
		0 & M\\ 
	\end{pmatrix}, ~~
	H_2= \begin{pmatrix}
		\sqrt{\nu} K& 0 \\ 
		0 & \sqrt{\nu} K\\ 
	\end{pmatrix},
\end{eqnarray*}
and
\begin{eqnarray*}
	S_{1}=\frac{1}{\theta}\begin{pmatrix}
		-i \omega \nu K & \sqrt{\nu} K\\ 
		-\sqrt{\nu} K & i \omega \nu K\\ 
	\end{pmatrix},~~ S_{2}=\begin{pmatrix}
		i \omega \sqrt{\nu} M & -M\\ 
		M & -i \omega \sqrt{\nu} M \\
	\end{pmatrix},
\end{eqnarray*}
\begin{eqnarray*}
	P_{1}=\frac{1}{\theta}\begin{pmatrix}
		I & -i \omega \sqrt{\nu} I\\ 
		i \omega \sqrt{\nu} I & -I\\ 
	\end{pmatrix},\quad
	P_{2}=\begin{pmatrix}
		0 & I\\ 
		I & 0\\ 
	\end{pmatrix},
\end{eqnarray*}
in which $\theta=1+\nu \omega^2$. In practice the authors suggest to use the matrix $V=H_1$ as a preconditioner and proved that the method is convergent under the condition $\alpha\geq \nu \omega^2/2$. The BAS iteration method induces the following preconditioner
\[
P_{BAS}(\alpha)=\zeta_{\alpha}
\begin{pmatrix}
	I &  (1+\omega^2\nu-i\omega\sqrt{\nu})I \\
	(1+\omega^2\nu+i\omega\sqrt{\nu})I  & -I
\end{pmatrix}
\begin{pmatrix}
	\alpha M+\sqrt{\nu} K  & 0 \\
	0  & \alpha M+\sqrt{\nu} K
\end{pmatrix}.
\] 
where $$\zeta_{\alpha}=\frac{(1+\alpha)}{\alpha(2+\omega^2\nu)}. $$ The reported numerical results in \cite{Zheng1} show that $\alpha=\theta$ and $\alpha=\theta/(1+\sqrt{\nu}\omega)$ are good choices for the BAS iteration and the BAS preconditioner, respectively.  

Recently, Salkuyeh in \cite{Salkuyeh} presented a stationary iteration method which was called Alternating Symmetric positive definite and Scaled symmetric positive semidefinite Splitting (ASSS) for solving the system (\ref{MAIN}). He used the idea of \cite{Zeng} and so rewrote the system \eqref{MAIN} as the 4-by-4 block real system
\begin{eqnarray*}
	{\mathcal{A}}  x \equiv\begin{pmatrix}
		M &  0 & \sqrt{\nu} K & \omega \sqrt{\nu} M\\
		0 & M & - \omega \sqrt{\nu} M & \sqrt{\nu} K\\
		\sqrt{\nu} K & - \omega \sqrt{\nu} M & -M & 0\\
		\omega \sqrt{\nu} M & \sqrt{\nu} K & 0 & -M
	\end{pmatrix} \begin{pmatrix}
		\mathfrak{R}(\bar{y}) \\
		\mathfrak{I}(\bar{y}) \\
		\mathfrak{R}(\bar{q}) \\
		\mathfrak{I}(\bar{q})
	\end{pmatrix} = \begin{pmatrix}
		\mathfrak{R}(\hat{y}_{d}) \\
		\mathfrak{I}(\hat{y}_{d}) \\
		0 \\
		0
	\end{pmatrix}\equiv\mathbf{\hat{b}}.
\end{eqnarray*}
The ASSS iteration method can be stated as
\begin{equation}\label{ASSS}
	\left\{
	\begin{array}{rl}
		(\alpha \mathcal{I}+{\mathcal M}){ x}^{(k+\frac{1}{2})}&=(\alpha \mathcal{I} -{\mathcal G} {\mathcal K}) { x}^{(k)}+b,\\
		(\alpha \mathcal{I}  +  {\mathcal K}){ x}^{(k+1)} &=(\alpha \mathcal{I} +\mathcal{G}\mathcal{M}) { x}^{(k+\frac{1}{2})}-\mathcal{G} b,
	\end{array}
	\right.
\end{equation}
where $\alpha>0$, $\mathcal{I}$ is an identity matrix of order $4m$, $\mathcal{K}=\sqrt{\nu} \mathcal{\hat{K}}/ \sqrt{\theta}$,
\begin{eqnarray*}
	\mathcal{M}=\begin{pmatrix}
		M &  0 & 0 & 0\\
		0 & M & 0 & 0\\
		0 & 0 & M & 0\\
		0 & 0 & 0 & M
	\end{pmatrix}, ~~ \mathcal{\hat{K}}=\begin{pmatrix}
		K &  0 & 0 & 0\\
		0 & K & 0 & 0\\
		0 & 0 & K & 0\\
		0 & 0 & 0 & K
	\end{pmatrix},
\end{eqnarray*}
and
\[
\mathcal{G}=\frac{1}{\sqrt{\nu \theta}} 
\left(\begin{array}{ccccc}
	0                   &  \omega \nu I    &  \sqrt{\nu} I      &  0 \\
	-\omega \nu I       &  0                 &  0   &  \sqrt{\nu} I  \\
	-\sqrt{\nu} I       & 0  &  0             &  -\omega \nu I    \\
	0   &  -\sqrt{\nu}I                 &  \omega \nu I                     &  0  \\ 
\end{array}\right).
\]
They proved that the ASSS iteration method is convergent unconditionally. Numerical results in \cite{Salkuyeh} show that $\alpha^{*}=\sqrt{\mu_{\min}\mu_{\max}}$ which $\mu_{\min}$ and $\mu_{\max}$ are the smallest and largest eigenvalues of the matrix $M$, respectively, is a good choice for the ASSS iteration method and the ASSS preconditioner. 

Throughout this paper, we use the following notations. For a given square matrix $A$,  the spectral radius and the spectrum of a $A$ are denoted by   $\rho(A)$ and $\sigma(A)$, respectively. We use  $\|.\|_2$ and $\|.\|_F$  for the Euclidean norm and the Frobenius norm, respectively.  The imaginary unit is denoted by $i=\sqrt{-1}$. The conjugate transpose of a matrix $A$, we use $A^H$.

In this paper, we present a new block alternating splitting  iteration  method (hereafter, we call it BASI‌ method) for solving the system (\ref{MAIN}) and prove that is convergent unconditionally. We also propose an estimation formula for the iteration parameter of the BASI method and the induced preconditioner.

\section{The BASI  method}\label{Sec2}
In this section, we propose a new splitting for the system (\ref{MAIN}). To do so, we define the matrices $\mathbf{S}_{1}$ and $\mathbf{S}_{2}$ as following
\begin{equation*}
	\mathbf{S}_{1}=
	\begin{pmatrix}
		I & -i \omega \sqrt{\nu} I \\ 
		i \omega \sqrt{\nu} I  & -I\\ 
	\end{pmatrix}, \quad \mathbf{S}_{2}=\begin{pmatrix}
		0 & \sqrt{\nu} I\\ 
		\sqrt{\nu}I & 0\\
	\end{pmatrix},
\end{equation*}
where $I \in \mathbb{R}^{m \times m}$ is the identity matrix. On the other hand, we have 
\begin{equation*}
	\mathbf{M}=
	\begin{pmatrix}
		M& 0 \\ 
		0 & M\\ 
	\end{pmatrix}, \quad \mathbf{K}=\begin{pmatrix}
		K & 0\\ 
		0 & K\\
	\end{pmatrix}.
\end{equation*}
Thus the system (\ref{MAIN}) can be written as following form 
\begin{equation}\label{EQ1}
	\mathbf{A}x \equiv \left( \mathbf{S}_{1} \mathbf{M} +\mathbf{S}_{2} \mathbf{K} \right)  x=b.
\end{equation}
We can see that 
$$\mathbf{S}_{1}^{H} \mathbf{S}_{1}=\theta I~~ \text{and}~~ \mathbf{S}_{2}^{H} \mathbf{S}_{2}=\nu I.$$ 
In addition, 
we can state the following lemma. 
\begin{lemma}
	If 
	\[
	\mathbf{S}=\frac{1}{\sqrt{\nu\theta}}\mathbf{S}_{1}^{H} \mathbf{S}_{2} ,
	\]
	then $S$ is skew-Hermitian. Moreover,    $\mathbf{S}^{2}=- \mathbf{I}$ where $\mathbf{I} \in \mathbb{R}^{2m \times 2m}$. 
	\begin{proof}
		By some computations, we see that 
		\begin{equation*}\label{Seq}
			\mathbf{S}=\frac{1}{\sqrt{\nu\theta}}
			\begin{pmatrix}
				-i \omega \nu I &  \sqrt{\nu} I\\ 
				-\sqrt{\nu} I & i \omega \nu I \\ 
			\end{pmatrix}.
		\end{equation*}
		Now, using $\theta=1+\nu \omega^2$  it is straightforward to show that   $\mathbf{S}^{2}=- \mathbf{I}$.
	\end{proof}
\end{lemma}

Multiplying both sides of (\ref{EQ1}) by $\mathbf{S}_{1}^{H}$, we obtain the system $$\left(  \mathbf{S}_{1}^{H} \mathbf{S}_{1} \mathbf{M}+\mathbf{S}_{1}^{H} \mathbf{S}_{2} \mathbf{K} \right) x = \mathbf{S}_{1}^{H}b,$$ which is equivalence to 
\begin{equation}\label{EQ2}
	\tilde{\mathbf{A}}x \equiv \left( \theta \mathbf{M} + \sqrt{\nu \theta}\mathcal{S} \mathbf{K}\right) x= \tilde{\mathbf{b}},
\end{equation}
where $\tilde{\mathbf{b}}=\mathbf{S}_{1}^{H} b$. For every $\alpha >0$, by adding $\alpha \mathbf{I}$ to the both sides of (\ref{EQ2}) we have
\begin{equation}\label{FIRST}
	\left( \alpha \mathbf{I} +\theta\mathbf{M}\right) x=\left( \alpha \mathbf{I}- \sqrt{\nu \theta} \mathbf{S}\mathbf{K}\right)  x+\tilde{\mathbf{b}}.
\end{equation}
On the other hand, by adding $\alpha \mathbf{S}$ to both sides of Eq. (\ref{EQ2}), we obtain the following equation
$$\left( \alpha \mathbf{S}+\sqrt{\nu \theta}\mathcal{S} \mathbf{K}\right) x=\left(  \alpha \mathbf{S} -\theta  \mathbf{M} \right)  x +\tilde{\mathbf{b}}.$$
Since $\mathbf{S}^{-1}=-\mathbf{S}$, by multiplying both sides of the above equation by $\mathbf{S}^{-1}$, we get
\begin{equation}\label{SECOND}
	\left( \alpha \mathbf{I}+\sqrt{\nu \theta} \mathbf{K} \right)  x=\left(  \alpha \mathbf{I} +\theta \mathbf{S} \mathbf{M}\right)  x-\mathbf{S} \tilde{\mathbf{b}}.
\end{equation}
Using Eqs. (\ref{FIRST}) and (\ref{SECOND}), we state the BASI  method for solving the system (\ref{EQ2}) as following.
\\ \\
\noindent\textbf{The BASI  method.}  {\it Let $ x^{(0)}\in \Bbb C ^{2n}  $ be an initial guess. For $k=0,1,2,\ldots $ until the sequence of iterates   $\{x^{(k)}\}_{k = 0}^\infty $ converges, compute the next iterate $ x^{(k+1)} $ via:
\begin{equation}\label{MBAS}
	\left\{ \begin{array}{l}
		\left( \alpha \mathbf{I} +\theta\mathbf{M}\right) x^{(k+\frac{1} {2})}=\left( \alpha \mathbf{I}- \sqrt{\nu \theta} \mathbf{S}\mathbf{K}\right)  x^{(k)}+\tilde{\mathbf{b}}, \\
		\left( \alpha \mathbf{I}+\sqrt{\nu \theta} \mathbf{K}\right) x^{(k+1)}=\left(  \alpha \mathbf{I} +\theta \mathbf{S} \mathbf{M} \right)  x^{(k+\frac{1}{2})}-\mathbf{S} \tilde{\mathbf{b}},
	\end{array}\right.
\end{equation}
where $ \alpha $ is a given positive constant. }
\bigskip

In each iteration of the BASI  method two systems with the coefficient matrices $ \alpha \mathbf{I} +\theta\mathbf{M}$ and $\alpha \mathbf{I}+\sqrt{\nu \theta} \mathbf{K}$ should be solved.  The system with the coefficient matrix $ \alpha \mathbf{I} +\theta\mathbf{M}$ can be split into two subsystems with the coefficient matrix $\alpha I +\theta M$, which is symmetric positive definite.  Hence, these subsystems  can be exactly solved  using the Cholesky factorization  or inexactly using the Conjugate Gradient (CG) method \cite{SaadBook}. Similarly, for solving the second half-step of the BASI method two subsystems with coefficient matrix $\alpha I+\sqrt{\nu \theta} K$  need to be solved, which can be treated  similar to the first half-step of the method.

By eliminating the vector $x^{(k+\frac{1}{2})}$ from Eq. (\ref{MBAS}) we get 
\begin{equation}\label{FIX}
	x^{(k+1)} = \mathbf{P}_{\alpha} x^{(k)}+\mathbf{Q}_{\alpha}\tilde{\mathbf{b}},
\end{equation}
where $$\mathbf{P}_{\alpha}=\left( \alpha \mathbf{I}+\sqrt{\nu \theta} \mathbf{K}\right) ^{-1} \left( \alpha \mathbf{I} +\theta \mathbf{S} \mathbf{M}\right)  \left( \alpha \mathbf{I} +\theta\mathbf{M}\right) ^{-1} \left( \alpha \mathbf{I}- \sqrt{\nu \theta} \mathbf{S}\mathbf{K}\right) ,$$
$$\mathbf{Q}_{\alpha}=\alpha \left( \alpha \mathbf{I}+\sqrt{\nu \theta} \mathbf{K}\right) ^{-1} \left( \mathbf{I}-\mathbf{S} \right)  \left( \alpha \mathbf{I} +\theta \mathbf{M}\right) ^{-1},$$
in which $\mathbf{P}_{\alpha}$ is the iteration matrix of  BASI  method. 
The next theorem investigates the convergence of the BASI  method.

\begin{theorem}\label{DKS}
	Let $M$, $K \in \mathbb{R}^{m \times m}$ be SPD matrices. Then for every $\alpha >0$, we have 
	$$\rho (\mathbf{P}_{\alpha}) \leq \eta_{\alpha}=\max_{\lambda \in \sigma(M)} \frac{\sqrt{\alpha^2+\theta^{2} \lambda^{2}
	}}{\alpha +\theta \lambda} \max_{\mu \in \sigma(K)} \frac{\sqrt{\alpha^{2}+\nu\theta \mu^{2}}}{\alpha+ \sqrt{\nu \theta} \mu} <1, $$
	where $\mathbf{P}_{\alpha}$ is the iteration matrix of the BASI method, $\nu $ and $\omega$ are given positive constants. This follows that the new iterative method converges unconditionally. 
	\begin{proof}
		It is easy to see that the matrix $\mathbf{P}_{\alpha}$ is similar to
		$$\tilde{\mathbf{P}}_{\alpha}= \left( \alpha \mathbf{I}+\sqrt{\nu\theta} \mathbf{K} \right)  \mathbf{P}_{\alpha} \left( \alpha \mathbf{I}+\sqrt{\nu\theta} \mathbf{K}\right) ^{-1}=\mathbf{U_{\alpha}} \mathbf{V_{\alpha}}$$ 
		where 
		$$\mathbf{U_{\alpha}}=\left( \alpha \mathbf{I} +\theta \mathbf{S}\mathbf{M}\right)  \left( \alpha \mathbf{I} +\theta \mathbf{M} \right) ^{-1}$$
		and
		$$\mathbf{V_{\alpha}}=\left( \alpha \mathbf{I} -\sqrt{\nu\theta} \mathbf{S} \mathbf{K}\right) \left( \alpha \mathbf{I}+\sqrt{\nu \theta}\mathbf{K}\right) ^{-1}.$$ 
		Then we have 
		\begin{eqnarray*}
			\rho(\mathbf{P}_{\alpha}) &\h=\h&\rho(\tilde{\mathbf{P}_{\alpha}}) \\
			&\h=\h& \rho(\mathbf{U_{\alpha}} \mathbf{V_{\alpha}}) \\
			&\h\leq\h& \parallel \mathbf{U_{\alpha}} \mathbf{V_{\alpha}} \parallel_{2} \\
			&\h\leq\h& \parallel \mathbf{U_{\alpha}} \parallel_{2} \parallel \mathbf{V_{\alpha}} \parallel_{2}.
		\end{eqnarray*} 
		On the other hand, we have 
		$\mathbf{S}^2=-\mathbf{I}$, $\mathbf{S}^{H}=-\mathbf{S}$, $\mathbf{S}\mathbf{M}=\mathbf{M}\mathbf{S}$ and $\mathbf{S}\mathbf{K}=\mathbf{K}\mathbf{S}$ 
		which follow
		\begin{eqnarray}
			\nonumber \| \mathbf{U_{\alpha}} \|^{2} _{2} &\h=\h&\rho\left( (\alpha \mathbf{I} +\theta \mathbf{M})^{-1} (\alpha \mathbf{I} -\theta \mathbf{M}\mathcal{S}) (\alpha \mathbf{I} + \theta \mathbf{S}\mathbf{M})(\alpha \mathbf{I} +\theta \mathbf{M})^{-1}\right) \\
			\nonumber &\h=\h& \rho \left( (\alpha \mathbf{I} +\theta \mathbf{M})^{-2} (\alpha^2 \mathbf{I} +\alpha \theta \mathbf{S}\mathbf{M}-\alpha \theta \mathbf{M}\mathbf{S}-\theta^{2} \mathbf{M} \mathbf{S}^{2} \mathbf{M})\right) \\
			&\h=\h& \rho \left( (\alpha \mathbf{I} +\theta \mathbf{M})^{-2} (\alpha \mathbf{I}+\theta^{2} \mathbf{M}^{2}) \right) \nonumber \\
			&\h=\h& \max_{\lambda \in \Bbb \sigma(M)} \frac{\alpha^2+\theta^{2} \lambda^{2}}{(\alpha +\theta \lambda)^{2}}. \label{Eqrho1}
		\end{eqnarray}
		We apply the same way for $\| \mathbf{V_{\alpha}} \|_{2}$ and get 
		\begin{eqnarray}
			\| \mathbf{V_{\alpha}} \|_{2} ^{2} &\h=\h& \max_{\mu \in \sigma(K)} \frac{\alpha^{2}+\nu\theta \mu^{2}}{\left( \alpha+ \sqrt{\nu \theta} \mu \right) ^{2}}. 
		\end{eqnarray}
		Since the matrices $M$ and $K$ are SPD, we have $\lambda$, $\mu>0$, for all $\lambda \in \sigma(M)$ and $\mu \in \sigma(K)$, accordingly we conclude that 
		$$\frac{\sqrt{\alpha^2+\zeta^{2} \lambda^{2}}}{\alpha +\zeta \lambda} <1,\quad  \text{and} \quad \frac{\sqrt{\alpha^{2}+\nu\theta \mu^{2}}}{\alpha+ \sqrt{\nu \theta} \mu}  <1.$$ Thus we can see that $ \| \mathbf{U_{\alpha}} \|_{2} <1$ and $ \| \mathbf{V_{\alpha}} \|_{2} <1$. So $$\rho (\mathbf{P}_{\alpha}) \leq \eta_{\alpha}=\| \mathbf{U_{\alpha}} \|_{2} \| \mathbf{V_{\alpha}} \|_{2} <1, $$ which compelets the proof.
	\end{proof}
\end{theorem}
Note that if the matrix $K$ is symmetric positive semidefinite (i.e., SPSD), then it is easy to see that the new iteration method is still convergence. 



\section{The BASI preconditioner}
It is not difficult to see that if  we define
$$\mathbf{B_{\alpha}} =\dfrac{1}{\alpha} \left( \mathbf{I}+\mathbf{S} \right) ^{-1} \left( \alpha \mathbf{I} +\theta\mathbf{M}\right)  \mathbf{S} \left( \alpha \mathbf{I}+\sqrt{\nu \theta} \mathbf{K}\right) ,$$
$$\mathbf{C_{\alpha}}  =  \dfrac{1}{\alpha} \left( \mathbf{I}+\mathbf{S}\right) ^{-1} \left(  \alpha \mathbf{S} -\theta \mathbf{M}\right)  \left( \alpha \mathbf{I}- \sqrt{\nu\theta} \mathbf{S}\mathbf{K} \right),$$
then 
$$\tilde{\mathbf{A}}=\mathbf{B_{\alpha}}-\mathbf{C_{\alpha}},\quad  \mathbf{P}_{\alpha}=\mathbf{B}_{\alpha}^{-1}\mathbf{C_{\alpha}}.$$
The latter equations show that the BASI  method induces the preconditioner $\mathbf{B}_{\alpha}$ for the system \eqref{EQ2}.
So if the BASI method is convergent, then the eigenvalues of  the matrix $\mathbf{B}_{\alpha}^{-1} \tilde{\mathbf{A}}$  are clustered in a circle with radius $1$, centred at $(1,0)$. So a Krylov subspace method like GMRES \cite{GMRES} or its flexible version (FGMRES) \cite{FGMRES} would be quite suitable for solving the preconditioned system (see \cite{BenziJCP})
\begin{equation}\label{Precon}
	\mathbf{B}_{\alpha}^{-1}\tilde{\mathbf{A}}x =\mathbf{B}_{\alpha}^{-1} \tilde{\mathbf{b}}.
\end{equation}
In each iteration of the GMRES method a vector of the form $w=\mathbf{B}_{\alpha}^{-1}v$ should be computed. Since, $\mathbf{S}^{-1}=-\mathbf{S}$,
we get
\[
w= -\alpha   ( \alpha \mathbf{I}+\sqrt{\nu \theta} \mathbf{K})^{-1}    \mathbf{S}  \left( \alpha \mathbf{I} +\theta\mathbf{M}\right)^{-1} \left( \mathbf{I}+\mathbf{S} \right) v.
\]
Hence, we can state the following algorithm for computing the vector $w$.\\

\smallskip
\noindent{\bf Algorithm 1:} Computation of $w=\mathbf{B}_{\alpha}^{-1}v$ \\
1. $p:=-\alpha\left( \mathbf{I}+\mathbf{S} \right) v.$\\
2. Solve $\left( \alpha \mathbf{I} +\theta\mathbf{M}\right)q=p$ for $q$.\\
3. $r:= \mathbf{S} q$.\\
4. Solve $( \alpha \mathbf{I}+\sqrt{\nu \theta} \mathbf{K})w=r$ for $w$.\\

In  the above algorithm, two systems with the coefficient matrices  $\alpha \mathbf{I} +\theta\mathbf{M}$ and $\alpha \mathbf{I}+\sqrt{\nu \theta} \mathbf{K}$. As we have discussed in the previous section, these systems can be solved exactly using the Cholesky factorization  or inexactly using the CG method. It is noted that if  these systems are solved inexactly, then the FGMRES method should be applied for the preconditioned system \eqref{Precon}.

\section{Parameter estimation}
In this section, using the idea of Ren and Cao \cite{Ren} we present a strategy for estimating the iteration parameter $\alpha$ of the new method preconditioner. We have  
$$\mathbf{C_{\alpha}}=\frac{1}{\alpha} \left( \mathbf{I}+\mathbf{S} \right)^{-1}  \left(  \alpha \mathbf{S} -\theta\mathbf{M} \right)  \left(\alpha \mathbf{I}- \sqrt{\nu \theta} \mathbf{S}\mathbf{K} \right) ,$$ and define the function $\varphi$ as 
$$\varphi(\alpha)=\frac{1}{\alpha} \| \left( \mathbf{I}+\mathbf{S} \right)^{-1}  \|_{F} \left( \alpha \| \mathbf{S} \|_{F} - \theta \|\mathbf{M}\|_{F}\right)  \left( \alpha \| \mathbf{I}\|_{F} -\sqrt{\nu \theta} \|\mathbf{S}\mathbf{K}\|_{F} \right).$$
To estimate the iteration parameter $\alpha$, we set 
$$\alpha \| \mathbf{S} \|_{F} - \theta \|\mathbf{M}\|_{F}=0.$$
Then, we obtain 
$$\alpha=\theta \frac{\|\mathbf{M}\|_{F}}{\|\mathbf{S}\|_{F}}.$$
Since 
$$\| \mathbf{S} \|_{F} ^{2}=trace(\mathbf{S}^{H} \mathbf{S})=trace(-\mathbf{S}^{2})=trace(\mathbf{I})=2m,$$
and 
$$\| \mathbf{M}\|_{F}^{2}=2 \|M\|_{F}^{2},$$ so the estimation formula for the iteration parameter $\alpha$ of the new method preconditioner is obtained as following 
\begin{equation}\label{est}
	\alpha_{est}=\theta \frac{\|\mathbf{M}\|_{F}}{\|\mathbf{S}\|_{F}}= \theta \frac{\sqrt{2}\|M\|_{F} }{\sqrt{2m}}=\theta \frac{\|M\|_{F} }{\sqrt{m}}.
\end{equation}
We will use (\ref{est}) in the numerical experiments and see that it gives often suitable results. Note that we can set $ \alpha \| \mathbf{I}\|_{F} -\sqrt{\nu \theta} \|\mathbf{S}\mathbf{K} \|_{F}=0$, but in practical implementation the formula (\ref{est}) is much better than that.
\section{Numerical results}
In this section, we use the numerical results to compare  the BAS, BASI and ASSS iteration methods for solving the complex linear system (\ref{MAIN}). For each method we present the number of iterations for the convergence and the elapsed CPU time (in seconds).  All of the numerical results are performed in \textsc{Matlab} R2018b by using a Laptop with 2.50 GHz central processing unit (Intel(R) Core(TM) i5-7200U), 6 GB RAM and Windows 10. In the Tables, we always use a zero vector as an initial guess and the maximum number of iterations is set to be 500 . A dagger ($\dagger$) means that the method has not convergence in 500 iterations.

We consider the distributed control problem introduced in Section \ref{Sec1} in two-dimensional case with the computational domain $\Omega=(0,1)\times (0,1)\in \Bbb{R}^2 $. The target state is set to be
\begin{equation}\label{Eq00108}
	y_d(x,y)=
	\left\{
	\begin{array}{cl}
		(2x-1)^2(2y-1)^2,&\text{if}~~(x,y)\in (0,\frac{1}{2})\times (0,\frac{1}{2}), \\
		0,&\text{otherwise}.
	\end{array} \right.
\end{equation}
In our numerical test, we discretize the problem using the bilinear quadrilateral {\bf Q1} finite elements with a uniform mesh \cite{Elman}.

In Tables \ref{Tbl1} and \ref{Tbl2}, first of all we compare the numerical results of the BASI  method described in Section \ref{Sec2} with the BAS and ASSS iteration methods. As we already mentioned,  in the implementation of the BASI  we need to solve  two subsystems with the coefficient matrix $\alpha I+ \theta M$ and the two subsystems with the coefficient matrix $\alpha I + \sqrt{\nu \theta} K$. In each iteration of the BAS method two subsystems with the coefficient matrix $\alpha M+\sqrt{\nu} K$ and two subsystems with the coefficient matrix $(1+\alpha ) M$ should be solved. In the ASSS iteration method we solve four subsystems with the coefficient matrix $\alpha I+M$ and four subsystems with the coefficient matrix $\alpha I +\sqrt{\frac{\nu}{\theta}} K$. All the systems are solved exactly using the sparse Cholesky factorization incorporated with the symmetric approximate minimum degree permutation. To do so, the \texttt{symamd.m} command of \textsc{Matlab} is applied. The outer iteration is stoped as soon as the residual norm of the system (\ref{MAIN}) is reduced by a factor of $10^{6}$.

After that the numerical results of the BASI preconditioner (denoted by P-BASI) described in Section \ref{Sec2} with the BAS and ASSS preconditioners (denoted by P-BAS and P-ASSS) described in Section \ref{Sec1} are compared. To do so, we use the complete version of the GMRES method in conjunction with aforementioned preconditioners.   All the subsystems are solved exactly using the sparse Cholesky factorization incorporated with the symmetric approximate minimum degree permutation. The iteration of  GMRES  as the outer iteration is stopped as soon as the residual norm is reduced by a factor of $10^6$.  We use the values of the estimation parameter $\alpha$ computed by the formula (\ref{est}) (presented in Tables \ref{Tbl01} and \ref{Tbl02}) for both the BASI method and BASI preconditioner.   Numerical results have been presented in Tables \ref{Tbl1} and \ref{Tbl2}.  In Tables \ref{Tbl1} and \ref{Tbl2}, we use $\alpha=\theta$ for the BAS method and apply $\alpha=\theta/(1+\sqrt{\nu}\omega)$ for the BAS preconditioner. We also use $\alpha^*$ (as suggested in \cite{Salkuyeh}) for both the ASSS method and the ASSS preconditioner.

\begin{table}[!t]
	\centering\caption{The values of an estimation and optimal parameters $\alpha$ of both the BASI  method and the BASI preconditioner for $h=2^{-7}$. \label{Tbl01}}
	\vspace{-0.4cm}
	\begin{center}
		\scalebox{0.77}{
			\begin{tabular}{|c|c|c|c|c|c|c|c|c|c|c|c|}
				\hline
				& \multicolumn{1}{c|}{$\nu \backslash \omega$} & \multicolumn{1}{c|}{$10^{-4}$} &  \multicolumn{1}{c|}{$10^{-3}$} & \multicolumn{1}{c|}{$10^{-2}$} &  \multicolumn{1}{c|}{$10^{-1}$} & \multicolumn{1}{c|}{$1$}& \multicolumn{1}{c|}{$10$}& \multicolumn{1}{c|}{$10^{2}$} &\multicolumn{1}{c|}{$10^{3}$}& \multicolumn{1}{c|}{$10^{4}$}	\\ \hline

				&$10^{-2}$& 0.00003 & 0.00003& 0.00003& 0.00003 & 0.000031& 0.000061&0.003080  &0.304939 &30.490909 \\
				$\alpha_{est}$&$10^{-4}$&0.00003 & 0.00003& 0.00003& 0.00003  &0.00003 &0.000031 &0.000061 & 0.003080 &0.304939\\
				&$10^{-6}$& 0.00003 & 0.00003 & 0.00003   &0.00003 & 0.00003 & 0.00003 &0.000031  &0.000061 &0.003080 \\
				&$10^{-8}$&0.00003 &0.00003 & 0.00003  & 0.00003 &0.00003 &0.00003  &0.00003 &0.000031 &0.000061\\ \hline 
				
				&$10^{-2}$& 0.00005 &  0.00005  & 0.00005&0.00005&0.00005  & 0.0001&0.005 & 0.4 &45  \\
				$\alpha_{opt}$&$10^{-4}$&0.00005 &0.00005 &0.00005 & 0.00005 & 0.00005 &0.00005 &0.00009 &0.004 &0.5\\
				&$10^{-6}$& 0.000045 & 0.000045 &  0.000045  & 0.000045&0.000045  & 0.000045 & 0.000045 &0.00009 &0.005 \\
				&$10^{-8}$& 0.00005&0.00005 & 0.00005  & 0.00005 &0.00005 & 0.00005 &0.00005 &0.00005 &0.0001\\ \hline
				
				&$10^{-2}$& 0.0001 &  0.0001  & 0.0001 &0.0001&0.0001  & 0.0001&0.005 & 0.4 &45  \\
				$\alpha^{*}_{opt}$&$10^{-4}$&0.00005 &0.00005 &0.00005 & 0.00005 & 0.00005 &0.00005 &0.00009 &0.004 &0.5\\
				&$10^{-6}$& 0.000045 & 0.000045 &  0.000045  & 0.000045&0.000045  & 0.000045 & 0.000045 &0.00009 &0.005 \\
				&$10^{-8}$& 0.00005&0.00005 & 0.00005  & 0.00005 &0.00005 & 0.00005 &0.00005 &0.00005 &0.0001\\ \hline
				
		\end{tabular}}
	\end{center}
\end{table}

\begin{table}[!t]
	\centering\caption{The values of an estimation and optimal parameters $\alpha$ of both the BASI  method and the BASI preconditioner for $h=2^{-6}$. \label{Tbl02}}
	\vspace{-0.4cm}
	\begin{center}
		\scalebox{0.85}{
			\begin{tabular}{|c|c|c|c|c|c|c|c|c|c|c|c|}
				\hline
				& \multicolumn{1}{c|}{$\nu \backslash \omega$} & \multicolumn{1}{c|}{$10^{-4}$} &  \multicolumn{1}{c|}{$10^{-3}$} & \multicolumn{1}{c|}{$10^{-2}$} &  \multicolumn{1}{c|}{$10^{-1}$} & \multicolumn{1}{c|}{$1$}& \multicolumn{1}{c|}{$10$}& \multicolumn{1}{c|}{$10^{2}$} &\multicolumn{1}{c|}{$10^{3}$}& \multicolumn{1}{c|}{$10^{4}$}	\\ \hline

				&$10^{-2}$&  0.00012 & 0.00012 & 0.00012 & 0.00012 & 0.00012 & 0.00024 & 0.01230  & 1.21867 & 121.8551 \\
				$\alpha_{est}$&$10^{-4}$& 0.00012& 0.00012 & 0.00012 & 0.00012  & 0.00012 & 0.00012 & 0.00024 & 0.01230  &1.21867 \\
				&$10^{-6}$& 0.00012 & 0.00012 & 0.00012   & 0.00012 & 0.00012 &0.00012 & 0.00012 & 0.00024 & 0.01230 \\
				&$10^{-8}$& 0.00012 & 0.00012 &   0.00012  &  0.00012 &0.00012 & 0.00012  & 0.00012 & 0.00012 &0.00024 \\ \hline 
				
				&$10^{-2}$& 0.0002 & 0.0002 & 0.0002 & 0.0002  & 0.0002 & 0.0004 & 0.02 & 1.7  & 123\\
				$\alpha_{opt}$&$10^{-4}$& 0.00015 & 0.00015 & 0.00015 & 0.00015  &0.00015 & 0.00015 & 0.0003 & 0.02  & 1.8\\
				&$10^{-6}$& 0.0002 & 0.0002 &  0.0002  & 0.0002 & 0.0002 & 0.0002& 0.0002 & 0.0004 & 0.025 \\
				&$10^{-8}$& 0.0002 & 0.0002 &  0.0002   & 0.0002  & 0.0002& 0.0002  & 0.0002 & 0.0002 & 0.0004\\ \hline

				&$10^{-2}$& 0.0002  & 0.0002 & 0.0002 & 0.0002 & 0.0002 & 0.0005 & 0.02  & 1.7 & 123 \\
				$\alpha^{*}_{opt}$&$10^{-4}$& 0.0002 & 0.0002  & 0.0002  &  0.0002  & 0.0002  & 0.0002  & 0.0004 & 0.02  & 1.7 \\
				&$10^{-6}$& 0.0002 & 0.0002 &  0.0002  & 0.0002 & 0.0002 & 0.0002 & 0.0002 & 0.0004 & 0.03 \\
				&$10^{-8}$& 0.0004 & 0.0004 &   0.0004  &  0.0004 & 0.0004 & 0.0004  & 0.0004 & 0.0004 &0.0009 \\ \hline

		\end{tabular}}
	\end{center}
\end{table}

\begin{table}[!t]
	\centering\caption{Number of iterations of the methods along with the elapsed CPU time (in parenthesis)  for $h=2^{-7}$ and different values of $\nu$ and $\omega$.\label{Tbl1}}
	\vspace{-0.4cm}
	\begin{center}
		\scalebox{0.75}{
			\begin{tabular}{|c|c|c|c|c|c|c|c|c|c|c|}
				\hline
				Method & \multicolumn{1}{c|}{$\nu \backslash \omega$} & \multicolumn{1}{c|}{$10^{-4}$} &  \multicolumn{1}{c|}{$10^{-3}$} & \multicolumn{1}{c|}{$10^{-2}$} &  \multicolumn{1}{c|}{$10^{-1}$} & \multicolumn{1}{c|}{$1$}& \multicolumn{1}{c|}{$10$}& \multicolumn{1}{c|}{$10^{2}$} &\multicolumn{1}{c|}{$10^{3}$}& \multicolumn{1}{c|}{$10^{4}$}	\\ \hline
				
				&$10^{-2}$&46(1.00)  &  46(1.00)   &46(1.00)  &46(0.99) &46(0.98)  &45(0.98) & 42(0.90) & 36(0.80) & 42(0.90)  \\
				BASI&$10^{-4}$& 42(0.96) &42(0.91) & 42(0.91) & 42(0.92) &42(0.90)  &42(0.90) &41(0.89) &36(0.79) &42(0.90)\\
				$(\alpha=\alpha_{est})$&$10^{-6}$& 36(0.83) & 36(0.80) &  36(0.80)  & 36(0.80) & 36(0.80)  & 36(0.80) & 36(0.79) & 37(0.82)&42(0.90) \\
				&$10^{-8}$&42(0.96) & 42(0.92)&  42(0.91) &42(0.91)  &42(0.91)  & 42(0.90) & 42(0.90)&42(0.90) &43(0.92)\\ \hline  
				
				&$10^{-2}$&40(0.98)  &  40(0.86)   & 40(0.86) &40(0.84) &40(0.86)  &39(0.88) & 37(0.85)& 34(0.84) &38(0.90)  \\
				BASI&$10^{-4}$&38(0.88) &38(0.84) & 38(0.84)&38(0.84)  & 38(0.83) &37(0.81) &37(0.86) &34(0.81) &38(0.90)\\
				$(\alpha=\alpha_{opt})$&$10^{-6}$& 34(0.80) & 34(0.77) & 34(0.77)   &34(0.75) &34(0.76)  & 34(0.76) &34(0.76)  & 34(0.82)&38(0.88) \\
				&$10^{-8}$&38(0.88) &38(0.84) & 38(0.83)  & 38(0.83) &38(0.84) & 38(0.87) &38(0.86) &  38(0.83)&38(0.88)\\ \hline

				&$10^{-2}$&  39(0.89) &  39(0.86)    &  39(0.85) & 39(0.87) &38(0.85) & 24(0.57) & 465(9.70) &$\dagger$ & $\dagger$    \\
				BAS&$10^{-4}$&  36(0.83)  & 36(0.82)  &  36(0.82)    & 36(0.80) & 36(0.80)  & 36(0.80)  & 39(0.87) & $\dagger$ & $\dagger$  \\
				&$10^{-6}$&  33(0.78)  & 33(0.75)  &  33(0.75)   & 33(0.76) &  33(0.75)  & 33(0.74) &  33(0.75) & 56(1.21) & $\dagger$  \\
				&$10^{-8}$& 38(0.89) & 38(0.86) &  38(0.85)  & 38(0.85) & 38(0.85)  & 38(0.85) &38(0.84) & 38(0.84)& 64(1.36)\\ \hline

				&$10^{-2}$&57(1.66)   &  57(1.53)  & 57(1.53)  & 57(1.55) & 57(1.54)& 56(1.51) & 53(1.45) & 44(1.20) &  51(1.38)   \\
				ASSS&$10^{-4}$& 53(1.46)   & 53(1.43)  &   53(1.45)   & 53(1.43) &  53(1.43) & 53(1.43)  &52(1.42) &44(1.21)  & 51(1.40)  \\
				&$10^{-6}$&  44(1.24)  & 44(1.22)  &  44(1.23)  &44(1.22) & 44(1.21) & 44(1.21) & 44(1.26)  & 43(1.19) & 51(1.38) \\
				&$10^{-8}$& 51(1.41) &51(1.39)  &  51(1.40)  & 51(1.38) &51(1.39) & 51(1.39) & 51(1.39)&51(1.39) &52(1.40) \\ \hline \hline \hline

				&$10^{-2}$& 31(1.60) & 31(1.50) & 31(1.48) &31(1.45) & 31(1.46) &31(1.44) &32(1.50) &34(1.56)& 28(1.34)  \\
				P-BASI&$10^{-4}$& 32(1.65) & 32(1.55) & 32(1.50) & 32(1.54) & 32(1.50) &32(1.50) &32(1.50) &34(1.58) &28(1.35) \\
				$(\alpha=\alpha_{est})$&$10^{-6}$& 32(1.62) & 32(1.56) &  32(1.50)  & 32(1.49) & 32(1.50)  &32(1.49) & 32(1.49) &32(1.49) & 28(1.35) \\
				&$10^{-8}$& 27(1.42)&27(1.35) & 27(1.33)  & 27(1.30) & 27(1.31) &27(1.30) &27(1.30) &27(1.30) &27(1.30)\\ \hline
				
				&$10^{-2}$& 25(1.40) &   25(1.39)  & 25(1.36) & 25(1.37) &25(1.40)&28(1.34) & 30(1.40) & 30(1.46)  & 24(1.38) \\
				P-BASI&$10^{-4}$&30(1.45) &30(1.39) &30(1.35) & 30(1.35) &30(1.35)  &30(1.36) &30(1.40) & 30(1.44)&24(1.30)\\
				$(\alpha=\alpha^{*}_{opt})$&$10^{-6}$&29(1.47)  & 29(1.43) & 29(1.40)   & 29(1.38)& 29(1.40) & 29(1.38) & 29(1.37) &30(1.48) & 24(1.30)\\
				&$10^{-8}$&23(1.28) &23(1.22) & 23(1.18)  & 23(1.19) &23(1.17) & 23(1.15) & 23(1.15)& 23(1.15)&22(1.24)\\ \hline  
				
				&$10^{-2}$& 18(0.94)  & 19(0.90)  & 20(0.93) & 20(0.91) &20(0.93)& 17(0.88) & 26(1.00) &49(1.73) & 46(1.86)  \\
				P-BAS&$10^{-4}$& 20(0.94)   & 21(0.90) &  22(0.92)   &22(0.92) & 22(0.94) & 22(0.89)  & 20(0.86) & 47(1.76)& 50(1.89) \\
				&$10^{-6}$& 18(0.89)  &  19(0.85) &  20(0.87)   &21(0.91)  &   21(0.91)  & 21(0.88) &  22( 0.91) &  28(1.07) &49( 1.82)  \\
				&$10^{-8}$& 18(0.90) & 19(0.84) & 20(0.87)  & 20(0.86) & 21(0.91)  &22(0.93) &22(0.91) &22(0.90) &28(1.08)\\ \hline 
				
				&$10^{-2}$& 36(1.43) & 36(1.39)  & 36(1.36) & 36(1.37) &36(1.42)  & 36( 1.35) & 38(1.43)&38(1.46) & 38(1.45)   \\
				P-ASSS&$10^{-4}$&  36(1.45)  & 36(1.39) &  36(1.35)   &36(1.35)  & 36(1.35) & 36(1.36) & 37(1.40)& 38(1.44) & 38(1.44) \\
				&$10^{-6}$& 37(1.47)  & 37(1.45) &   37(1.40)  & 37(1.39) & 37(1.40)  & 37(1.40) & 37(1.38) & 38(1.48) & 38(1.43)\\
				&$10^{-8}$& 37(1.47) & 37(1.43) & 37(1.40)  & 37(1.40) & 37(1.40)  &  37(1.40) &37(1.40) &37(1.39) &36(1.36)\\ \hline

		\end{tabular}}
	\end{center}
\end{table}
\begin{table}[!t]
	\centering\caption{Number of iterations of the methods along with the elapsed CPU time (in parenthesis)  for $h=2^{-6}$ and different values of $\nu$ and $\omega$.\label{Tbl2}}
	\vspace{-0.4cm}
	\begin{center}
		\scalebox{0.7}{
			\begin{tabular}{|c|c|c|c|c|c|c|c|c|c|c|}
				\hline
				Method & \multicolumn{1}{c|}{$\nu \backslash \omega$} & \multicolumn{1}{c|}{$10^{-4}$} &  \multicolumn{1}{c|}{$10^{-3}$} & \multicolumn{1}{c|}{$10^{-2}$} &  \multicolumn{1}{c|}{$10^{-1}$} & \multicolumn{1}{c|}{$1$}& \multicolumn{1}{c|}{$10$}& \multicolumn{1}{c|}{$10^{2}$} &\multicolumn{1}{c|}{$10^{3}$}& \multicolumn{1}{c|}{$10^{4}$}	\\ \hline
				
				&$10^{-2}$& 45(0.25)  &   45(0.19)   & 45(0.20)  & 45(0.20) & 45(0.20) & 44(0.21) & 40(0.19) & 35(0.16) & 43(0.19)  \\
				BASI&$10^{-4}$& 40(0.25) &40(0.19) & 40(0.19) & 40(0.18)  &  40(0.18) & 40(0.18) &39(0.18) & 35(0.16) & 43(0.20) \\
				$(\alpha=\alpha_{est})$&$10^{-6}$& 35(0.22) & 35(0.17)  & 35(0.17) & 35(0.17) & 35(0.16)  & 35(0.16) & 35(0.17)  & 36(0.17) & 43(0.18)  \\
				&$10^{-8}$& 43(0.24) & 43(0.18) &  43(0.18) & 43(0.19) & 43(0.20)  & 43(0.19) & 43(0.18) & 43(0.19) & 43(0.19)\\ \hline  
				
				&$10^{-2}$&  40(0.22)  & 40(0.20) & 40(0.19) & 40(0.19) & 40(0.17) & 40(0.17) & 38(0.16) & 33(0.16)  & 42(0.19) \\
				BASI&$10^{-4}$& 37(0.17) &37(0.16) & 37(0.18) & 37(0.16)  & 37(0.18)  & 37(0.18) & 36(0.17) & 33(0.16) & 38(0.18)\\
				$(\alpha=\alpha_{opt})$&$10^{-6}$& 33(0.20) & 33(0.18) & 33(0.16) & 33(0.17) & 33(0.18)  & 33(0.16) & 33(0.17)  & 33(0.16) &38(0.17)  \\
				&$10^{-8}$& 38(0.23) & 38(0.18) & 38(0.17)  & 38(0.17) &  38(0.17) & 38(0.16) &38(0.17) & 38(0.18) & 38(0.17)\\ \hline

				&$10^{-2}$& 38(0.30)   &  38(0.20)   &  38(0.20)  & 38(0.21) & 38(0.20)  & 24(0.16)  & 476(1.82)  & $\dagger$  & $\dagger$      \\
				BAS&$10^{-4}$& 35(0.30)  &35(0.20) & 35(0.21) & 35(0.21) & 35(0.20) & 35(0.20)   & 39(0.23)  & $\dagger$  & $\dagger$    \\
				&$10^{-6}$& 33(0.29)   & 33(0.20) & 33(0.21)  & 33(0.23)  & 33(0.21) & 33(0.19) & 56(0.30)  &  $\dagger$ & $\dagger$ \\
				&$10^{-8}$&38(0.30)  & 38(0.20) &  38(0.20)  & 38(0.20) &38(0.20) &38(0.20) &38(0.20) & 38(0.19) & 65(0.31) \\ \hline

				&$10^{-2}$& 56(0.35)  & 56(0.33) &  56(0.32)  & 56(0.30) & 56(0.30) & 55(0.31)  & 50(0.28) &40(0.23) & 51(0.28)    \\
				ASSS&$10^{-4}$&  50(0.31)   &50(0.29) &  50(0.28)  &50(0.28) & 50(0.30)  & 50(0.28) & 48(0.28) & 40(0.25)  &  51(0.29)  \\
				&$10^{-6}$&  48(0.28) & 40(0.23)& 40(0.23)  & 40(0.25) & 40(0.23) & 40(0.23) &  40(0.25)  & 42(0.24) & 51(0.29)  \\
				&$10^{-8}$& 51(0.33)  & 51(0.28) &  51(0.29)  &51(0.28) & 51(0.28) & 51(0.30)  &51(0.28) &51(0.27) &52(0.28)  \\ \hline \hline \hline

				&$10^{-2}$&  32(0.35) & 32(0.34)  & 32(0.30) &32(0.30) &  32(0.29)& 32(0.30) &34(0.31) &32(0.28) & 26(0.25)   \\
				P-BASI&$10^{-4}$& 32(0.39) &  32(0.32)  &  32(0.32) & 32(0.29)  &  32(0.31) & 32(0.29) & 32(0.32) & 32(0.30) & 26(0.23) \\
				$(\alpha=\alpha_{est})$&$10^{-6}$&  31(0.38) &  31(0.30) &  31(0.30)  &  31(0.29) &  31(0.29)  & 31(0.29) &  31(0.29) & 30(0.28) & 26(0.25) \\
				&$10^{-8}$& 24(0.32) & 24(0.24) &  24(0.23) &  24(0.24)  &  24(0.23) & 24(0.22) & 24(0.23) & 24(0.23) &24(0.22) \\ \hline

				&$10^{-2}$&  25(0.34) & 25(0.24)  & 25(0.25) & 25(0.25) & 25(0.22) & 28(0.28) & 30(0.30) & 28(0.31) &  24(0.29)  \\
				P-BASI&$10^{-4}$& 30(0.35) &  30(0.30)  & 30(0.31)  & 30(0.28)  & 30(0.29)  & 30(0.27) & 30(0.28) & 28(0.32) & 20(0.30) \\
				$(\alpha=\alpha_{opt})$&$10^{-6}$& 27(0.34)  &  27(0.27) &  27(0.27)  & 27(0.25)  &  27(0.27)  & 27(0.24) &  27(0.26) & 27(0.30) & 16(0.24) \\
				&$10^{-8}$& 16(0.25) & 16(0.18) &  16(0.17) & 16(0.17)   & 16(0.18)  & 16(0.17) & 16(0.17) & 16(0.15) & 14(0.15)\\ \hline

				&$10^{-2}$& 18(0.30)  & 19(0.20)  & 20(0.20) &20(0.20) & 20(0.21) & 17(0.18)  & 26(0.25) & 49(0.50) &  44(0.40)  \\
				P-BAS&$10^{-4}$& 20(0.31) &  20(0.22)  & 21(0.25)  & 22(0.21)  &  22(0.25) & 22(0.23) & 20(0.21) & 47(0.45) & 48(0.45) \\
				&$10^{-6}$& 18(0.30)  & 19(0.23)  &  20(0.22)  & 21(0.22)  &   21(0.22) & 21(0.23) &  22(0.23) & 28(0.30) & 47(0.50) \\
				&$10^{-8}$& 17(0.30) & 19(0.22) &  20(0.22) &  20(0.24)  & 21(0.25)  &21(0.22) & 22(0.25) & 22(0.23) & 28(0.24) \\ \hline
				
				&$10^{-2}$& 37(0.42)  & 37(0.40)  & 37(0.35) &37(0.33) & 37(0.36) & 38(0.35) & 38(0.35)& 36(0.34) &  32(0.29)  \\
				P-ASSS&$10^{-4}$& 37(0.41) &  37(0.38)  &  37(0.34) & 37(0.34)  & 37(0.35)  & 37(0.34) & 38(0.36) & 38(0.35) & 32(0.31) \\
				&$10^{-6}$&  35(0.42) & 35(0.35)  &  35(0.32)  &  35(0.32) &  35(0.34)  & 35(0.31) &  35(0.33) & 36(0.34) & 32(0.30) \\
				&$10^{-8}$& 31(0.37) & 31(0.30) & 31(0.31)  & 31(0.29)   & 31(0.30)  & 31(0.30)& 31(0.30) & 31(0.30) & 30(0.29)\\ \hline

		\end{tabular}}
	\end{center}
\end{table}

\section{Conclusion}\label{Sec5}
We have presented a modified  version of the block alternating splitting iteration (BASI)  method  for solving the system arising from finite element discretization of the distributed optimal control problem with time-periodic parabolic equations. We have proved that the proposed method is unconditionally convergent. An  estimation formula for the iteration parameter of the BASI preconditioner has been given. Numerical results show that  both the BASI  method and the BASI preconditioner are efficient and outperform the BAS iteration method.


\end{document}